\documentstyle[psfig]{article}

\def\0{{\bf 0}}

\def\bt{\begin{theorem}}
\def\et{\end{theorem}}
\def\bp{\begin{proposition}}
\def\ep{\end{proposition}}
\def\bl{\begin{lemma}}
\def\el{\end{lemma}}
\def\bi{\begin{itemize}}
\def\ei{\end{itemize}}
\def\bd{\begin{description}}
\def\ed{\end{description}}
\def\br{\begin{remark}}
\def\er{\end{remark}}
\def\be{\begin{equation}}
\def\ee{\end{equation}}
\def\bc{\begin{corollary}}
\def\ec{\end{corollary}}
\def\bex{\begin{example}}
\def\eex{\end{example}}

\input mssymb12.tex

\newcommand{\Z}{{\Bbb Z}}
\newcommand{\R}{{\Bbb R}}
\newcommand{\C}{{\Bbb C}}

\def\NABLA#1{{\mathop{\nabla\kern-.5ex\lower1ex\hbox{$#1$}}}}
\def\Nabla#1{\nabla\kern-.5ex{}_#1}

\newtheorem{theorem}{Theorem} 
\newtheorem{corollary}[theorem]{Corollary}
\newtheorem{lemma}[theorem]{Lemma}
\newtheorem{proposition}[theorem]{Proposition}

\newtheorem{remark}[theorem]{Remark}
\newtheorem{example}[theorem]{Example}

        \newfont{\Frak}{eufm10}                 
        \newcommand{\frak}[1]{\mbox{\Frak {#1}}}

\title{The Geometry of Grauert Tubes and Complexification of Symmetric Spaces}

\author{D. Burns\thanks{Partially supported by NSF, DMS-0104047.} \and
S. Halverscheid \and R. Hind}

\begin{document}

\maketitle

\vskip 5mm

\noindent{\bf{Abstract}}

\vskip 3mm

\noindent We study the canonical complexifications of non-compact
Riemannian symmetric spaces by the Grauert tube construction. We
determine the maximal such complexification, a domain already
constructed by Akhiezer and Gindikin \cite{ag}, and show that this
domain is Stein. We also determine when invariant complexifications,
including the maximal one, are Hermitian symmetric. This is expressed
simply in terms of the ranks of the symmetric spaces involved.

\vskip 5mm

\section{Introduction}

The purpose of this paper is to investigate certain canonical
complexifications of irreducible Riemannian symmetric spaces of
the noncompact type. We will be interested in determining whether
these complexifications are rigid in the sense that their complex
automorphism group is isomorphic to the isometry group of the
symmetric space (this is always a subgroup) or whether the
complex manifold exhibits additional symmetries. We will also deal
with a question of Steinness, verifying a conjecture of D. Akhiezer
and S. Gindikin.

The construction of the canonical complexifications is actually
very general and applies to any real-analytic Riemannian manifold.
A general rigidity result was proven in our previous paper
\cite{bh} for the complex manifolds associated to compact
Riemannian manifolds, but almost nothing is known about the
complexifications of arbitrary noncompact manifolds. In the
case of symmetric spaces we are continuing work of S-J. Kan
and D. Ma from \cite{km}.

We start with the following theorem proved independently
by V. Guillemin and M. Stenzel in \cite{gs} and by L. Lempert and R. Sz\"{o}ke in \cite{ls}.

Suppose that $(M,g)$ is a real-analytic Riemannian manifold of dimension $n$. Identify $M$ with the zero section in $TM$. 
Let
$\rho :TM\to \R$ be the length, with respect to $g$, of tangent
vectors.

\begin{theorem}  There exists a unique complex structure on all sufficiently
small neighbourhoods of $M$ in $TM$ such that the following conditions hold.

$(i)$ $\rho ^2$ is strictly plurisubharmonic and the corresponding K\"{a}hler
metric restricts to $g$ on $M$.

$(ii)$ $\rho$ is a solution of the homogeneous complex Monge-Amp\`{e}re equation $(dd^c\rho)^n = 0$ on
$T^r M\setminus M$, where $M \subset T^r M$.
\end{theorem}

Let $T^r M=\{v\in TM | \rho (v) < r\}$. Then when $M$ is a compact
manifold or a symmetric space the complex structure exists on the
disk bundles $T^r M$ for $r$ sufficiently small. In this case we
call the resulting complex manifolds $T^r M$ Grauert tubes, and
$r$ will often be called the radius of the tube.
We will
also call the complex manifolds Grauert domains even if the underlying
subset of $T^r M$ is not necessarily a tube. For a general noncompact
$M$ the canonical complex structure need not exist on any $T^r M$.

An equivalent characterization of these complex manifolds (see also
\cite{ls}) is that the complex structure is the unique one making
the leaves of the Riemann foliation with their natural complex
structures into holomorphic curves. In other words, for any geodesic
$\gamma :\Bbb R \to M$ the map
$$d\gamma : \Bbb C \to TM$$
$$s+it \mapsto (\gamma (s), t\gamma '(s))$$
is actually holomorphic with respect to the Grauert domain structure,
wherever that structure is defined.

We remark here that all Grauert tubes have an antiholomorphic involution
$\sigma:T^r M\to T^r M$ given by $v\mapsto -v$.

By functoriality of the construction, the differential of any
isometry of $M$ is actually a biholomorphism of a Grauert tube
$T^r M$. We will say that a Grauert domain $N$, invariant under
the action of $Isom(M,g)$, is {\it rigid} if its complex
automorphism group $Aut_{\Bbb C}(X)=dIsom(M,g)$.

The following theorem is taken from \cite{bh}.

\begin{theorem}
Any Grauert tube $T^r M$ of finite radius $r$
over a compact manifold $M$ is rigid.
\end{theorem} 

Now, for symmetric (or locally symmetric) spaces of nonnegative
curvature the canonical complex structure exists on the whole
tangent bundle $TM$ (see for instance \cite{sz}). The resulting
complex manifolds for compact rank one spaces are described by
G. Patrizio and P-M. Wong in \cite{pmw}. As a particular case, the Grauert
tube complex structure on $TS^2$ associated to the standard
round metric on $S^2$ is biholomorphic to the complex quadric
$\{z_1^2+z_2^2+z_3^2=1\}$ in $\Bbb C^3$. Thus such tubes of
infinite radius are not rigid in general.

For Riemannian manifolds with some negative sectional curvature,
the whole tangent bundle cannot be given the structure of a
Grauert domain. In fact the following theorem of Lempert and
Sz\"{o}ke (see \cite{ls} again) gives an upper bound on the radius
of a Grauert tube.

\begin{theorem}
Suppose that the sectional curvature of a $2$-plane of
$(M,g)$ is equal to $-\lambda$, with $\lambda >0$.
Then if a Grauert tube
structure exists on $T^r M$ the radius satisfies
$$r\le \frac{\pi}{2\sqrt{-\lambda}}.$$  
\end{theorem}

Hence for manifolds $(M,g)$ with some negative sectional curvature
there is a finite radius $r_{max}$ for which a Grauert tube
structure exists on $T^r M$ if and only if $r\le r_{max}$.
We emphasize that Theorem $2$ holds even if $r=r_{max}$.

From now on we will concentrate on the case when $(M,g)$ is
a symmetric space of the noncompact type.
The following theorem is proven by Kan and Ma in \cite{km}.

\begin{theorem}
Let $(M,g)$ be a symmetric space and $T^r M$ an associated
Grauert tube with $r<r_{max}$. Then $T^r M$ is either rigid
or the ball.
\end{theorem}

It turns out that there is at most one $r$ for which $T^r M$
is covered by a ball, so we can say that all but at most one
$T^r M$ are rigid for $r<r_{max}$. One surprising result of
this paper is that for rank one symmetric spaces
of noncompact type, the tubes $T^{r_{max}} M$ are never rigid,
in fact they are all Hermitian symmetric.

We define the {\it maximal Grauert domain} associated to a
real-analytic $(M,g)$ to be the largest connected domain
in $TM$ containing $M$ on which our canonical complex structure
exists, that is, on which we can define a complex structure
satisfying the conditions of Theorem $1$.

Without loss of generality suppose that $(M,g)$ is an irreducible
symmetric space of noncompact type. It is not hard to show
that when $(M,g)$ is of rank one, the maximal Grauert domain
is just the Grauert tube of maximal radius. However for
higher rank cases the maximal Grauert domain is always
larger.

In section $2$ we will show that the maximal Grauert domains
can be described algebraically, and correspond to the domains
defined and studied by D. Akhiezer and S. Gindikin amongst others,
see \cite{ag}.

Let $M$ be written in Klein form as $M=G/K$ where $G$ is a
connected semisimple Lie group with finite center and $K$
is a maximal compact subgroup. The groups $G$ and $K$ can
be complexified to linear algebraic groups $G^{\Bbb C}$ and
$K^{\Bbb C}$ over $\Bbb C$. Let $o=e.K^{\Bbb C}\in G^{\Bbb C}/K^{\Bbb C}$.

Now let $\frak g$ be the Lie algebra of $G$ and
$\frak g= \frak t + \frak p$ be a Cartan decomposition corresponding
to the pair $(G,K)$. We let $\frak a$ be a maximal abelian
subalgebra of $\frak p$ and $\Sigma$ be a root system of
$\frak g$ with respect to $\frak a$.

We will identify $G(o)$ with $M=G/K\hookrightarrow G^{\Bbb C}/K^{\Bbb C}$.
Each $G$-orbit in $G^{\Bbb C}/K^{\Bbb C}$ intersects
$exp(i\frak a).o$ in an orbit of the Weyl group.

Finally, following \cite{ag} we define
$$\omega=\{H\in\frak a | |\alpha(H)|<\frac{\pi}{2} \forall \alpha \in \Sigma\}.$$

\begin{theorem}
The maximal Grauert domain associated to $M=G/K$ is biholomorphic
to the domain $D=G(exp(i\omega)).o\subset G^{\Bbb C}/K^{\Bbb C}$.
\end{theorem}

This is essentially contained in \cite{ag} or \cite{sz}
but we clarify it in section $2$.

There is an alternative for the complex automorphism group
of maximal Grauert domains. We prove that the domains are
either rigid or biholomorphic to Hermitian symmetric spaces.
In fact we will prove the following.

\begin{theorem}
Let $\omega' \subseteq \omega$ be a symmetric convex domain
invariant under the action of the Weyl group. Then the
domain $D'=G(exp(i\omega')).o$ is either Hermitian symmetric
or $Aut_{\Bbb C}(D')\cong Isom(M)$.
\end{theorem}

It remains to decide in specific cases whether a maximal Grauert
domain is rigid or Hermitian symmetric. In the second case,
we will have a Hermitian symmetric space with an antiholomorphic
involution whose fixed point set is isomorphic to the original
symmetric space. Such objects have been classified by H. Jaffee,
see \cite{jaf1} and \cite{jaf2}. An immediate consequence is
that the maximal Grauert domain of $SL(3,\Bbb R)/SO(3)$ is rigid.

Suppose then that we have a Hermitian symmetric space $N$ with
an anti-holomorphic involution $\sigma$ whose fixed point set
is our symmetric space $M=G/K$. Then $N$ is always a Grauert
domain corresponding to $M$. We need to decide whether $N$ is
the maximal domain. Suppose $N=G'/K'$ in Klein form and
$\frak g'=\frak t'+\frak p'$ is its Cartan decomposition.
Again letting $\frak g=\frak t+\frak p$ be the decomposition
of the Lie algebra of $G$, since $M\subset N$ is totally geodesic
we have $\frak t \subset \frak t'$ and $\frak p \subset \frak p'$.
Let $\frak a$ be a maximal abelian subalgebra of $\frak p$ and
$\frak a'$ with $\frak a \subset \frak a' \subset \frak p'$ be
a maximal abelian subalgebra of $\frak p'$. 
Let $\Sigma \subset \frak a^*$ be the root system of $\frak g$ with
respect to $\frak a$ and
$\Sigma' \subset \frak a'^*$ be the root system of $\frak g'$ with
respect to $\frak a'$.

\begin{theorem}
The Hermitian space $N$ is the maximal Grauert domain corresponding
to $M$ if and only if
$$\max_{\alpha\in\Sigma}|\alpha (H)|=\max_{\alpha\in\Sigma'}|\alpha (H)|\forall H\in\frak a.$$     
\end{theorem}

For example, if $M$ itself is Hermitian symmetric, then the maximal
Grauert domain is biholomorphic to the product $M\times \overline{M}$
into which $M$ embeds diagonally. A result close to this special case
of the theorem is contained in \cite{wz}.

The condition in Theorem $8$ is sometimes awkward to check, especially
in exceptional cases, but fortunately it is equivalent to the simple
condition given in the next proposition.

\begin{proposition}
The Hermitian space $N$ is the maximal Grauert domain corresponding
to $M$ if and only if $rank(N)=2rank(M)$.
\end{proposition}

Combining Proposition $8$ with H. Jaffee's classification, it is easy to obtain
a general description of the maximal Grauert domains of irreducible
symmetric spaces of the noncompact type. This section will conclude
with the resulting list.

Whether they are Hermitian symmetric or not, we are able to prove
the following general result on the structure of the maximal
Grauert domains.

\begin{theorem}
The maximal Grauert domain associated to a symmetric space is Stein.
\end{theorem}

This result was conjectured by D. Akhiezer and S. Gindikin in
\cite{ag}, who showed some examples of this property. These domains
are intimately related to the linear cycle spaces of P.A. Griffiths,
see \cite{bilbao} for a review. Several examples of the theorem were
thus shown by Huckleberry, Wolf and Zireau in various papers on linear
cycle spaces (see \cite{bilbao} for references), using the fact that
the linear cycle spaces are Stein, see \cite{w1}, for example. In
particular, the results of the current paper were in part motivated by
the appearance of ``hidden symmetries'' for some linear cycle spaces,
as pointed out to us by Joe Wolf.

The proof of the last theorem is by a direct characterization of the
$G$-invariant plurisubharmonic functions on the maximal Grauert
domain. Earlier examples of this result are due to Azad and Loeb
\cite{azad} for compact symmetric spaces, and K.-H. Neeb \cite{neeb}
for certian non-degenerate semigroups.

After giving the algebraic description of maximal Grauert domains in
section $2$, we show in section $3$ that for a Hermitian symmetric space
the corresponding maximal domain is simply the product of the space
with itself. This is used to prove Theorem $7$ in section $4$. 
In section $5$ the maximality condition is shown to be equivalent
to the simpler statement in Proposition $8$.
In section $6$ we prove the alternative for the biholomorphism
groups of the Grauert domains. Finally in section $7$ we prove that
the maximal Grauert domains are Stein by characterizing, as already
noted, their $G$-invariant plurisubharmonic functions.  

The third author would like to thank Brian Hall for several
enlightening discussions on these topics.

\pagebreak

{\bf Maximal Grauert domains of noncompact symmetric spaces}

\medskip

\begin{tabular}{rl}

$SL(n,\Bbb R)/SO(n), N>2$ & rigid \\
$SU^*(2n)/Sp(n)$ & rigid \\
$SU(p,q)/S(U_p \times U_q)$ & product \\
$SO_0(2,1)/SO(2)$ & product \\
$SO_0(p,1)/SO(p), p>2$ & $SO_0(p,2)/SO(p)\times SO(2)$ \\
$SO_0(p,2)/SO(p)\times SO(2)$ & product \\
$SO_0(p,q)/SO(p)\times SO(q), q>2$ & rigid \\
$SO^*(2n)/U(n)$ & product \\
$Sp(n,\Bbb R)/U(n)$ & product \\
$Sp(p,q)/Sp(p)\times Sp(q)$ & $SU(2p,2q)/S(U_{2p}\times U_{2q})$ \\
$SL(n,\Bbb C)/SU(n), n>2$ & rigid \\
$SO(n,\Bbb C)/SO(n), n>3$ & rigid \\
$Sp(n,\Bbb C)/Sp(n), n>1$ & rigid \\
$(\frak{e}_{6(-14)}, \frak{so}(10)+\Bbb R)$ & product \\
$(\frak{e}_{7(-25)}, \frak{e}_6 +\Bbb R)$ & product \\
$(\frak{f}_{4(-20)}, \frak{so}(9))$ & $(\frak{e}_{6(-14)}, \frak{so}(10)+\Bbb R)$ \\
all other exceptional spaces & rigid \\

\end{tabular}

\medskip

The notation above is taken from the book \cite{helg}. We notice that
all rank $1$ examples have Hermitian symmetric maximal Grauert tubes.
Some real symmetric spaces do appear as the fixed point sets of
involutions on Hermitian symmetric spaces but nevertheless have
rigid maximal domains. Examples of this are the spaces
$SO_0(p,q)/SO(p)\times SO(q)$, for $q>2$ which appear as the fixed
point sets of involutions on $SU(p,q)/S(U_p \times U_q)$. 
However, if $p$ and $q$ are even, then $SU(p,q)/S(U_p \times U_q)$
is a maximal Grauert domain. Another similar example is
$SU^*(8)/Sp(4)$ inside $(\frak{e}_{7(-25)}, \frak{e}_6 +\Bbb R)$.

\section{Algebraic description of maximal Grauert domains}

For a noncompact symmetric space $(M,g)$, let
$G$, $K$, $G^{\Bbb C}$, $K^{\Bbb C}$, $\frak g$, $\frak t$,
$\frak p$, $\frak a$ be as in the introduction, that is, 
$M=G/K$ where $G$ has Lie algebra $\frak g=\frak t+\frak p$ and
$\frak a \subset \frak p$ is maximal abelian.
Let $\Sigma$ be a root system of $\frak g$ with respect to $\frak a$
for the decomposition
\begin{equation}
\label{root-decomposition}
\frak g = Z(\frak a) \oplus \bigoplus_{\alpha \in \Sigma} g^{\alpha}
\end{equation}
of the adjoint representation of $\frak a$.
Here, $Z(\frak a)$ denotes the center of $\frak a$.

We can identify the tangent bundle $TM$ with $G\times_K \frak p$ where
for $k\in K$, $g\in G$ and $X\in \frak p$,
$$k(g,X)=(gk^{-1},Ad(k)X).$$

With this identification $G\times_K (Ad(K)\omega)$ is a subdomain of
$TM$ containing $M$, where
$$
\omega=\{H\in\frak a | |\alpha(H)|<\frac{\pi}{2} 
                            \forall \alpha \in \Sigma\}.
$$
As in the introduction, let $D=G.exp(i\omega).o\subset G^{\Bbb C}/K^{\Bbb C}$.

Then it is proven in \cite{ag} that the map
$$\phi:G\times_K (Ad(K)\omega)\to D$$
$$(g,Ad(k)H)\mapsto gk \exp(iH).o$$
is a real-analytic $G$-equivariant diffeomorphism.

Furthermore, the leaves of the Riemann foliation map onto
holomorphic curves in $D$ and so the pull-back of the complex
structure on $D$ gives $G\times_K (Ad(K)\omega)$ the (unique)
adapted complex structure.

It remains to show that this is actually the maximal domain.
To do this, suppose that $H\in \omega$ and that there exists
some $\alpha \in \Sigma$ with $\alpha (H)=\frac{\pi}{2}$.
Let $\gamma$ be the geodesic in $M=G/K$ with $\gamma'(0)=H$.
Then we have a map
$$d\gamma:\Bbb C\to TM$$
$$s+it\mapsto (\gamma(s),t\gamma'(s)).$$

We need to show that a Grauert domain structure can extend only
over $d\gamma (|t|<1)$. In \cite{sz}, section $3$, a precise
criterion describes how far the complex structure can be
extended for symmetric spaces. In the following, this criterion
will be related to the root decomposition (\ref{root-decomposition}).

Note that the linear operator on $\frak p$ 
$$ Y \mapsto R(Y,H)H := - (\mathrm{ad} H)^2(Y)$$
is the Jacobi operator for the geodesic 
$\gamma$ in $T_{eK} G/K \cong \frak p$. 
Since $Z(\frak a) \cap \frak p = \frak a$,
the decomposition into eigenspaces of the Jacobi operator
is given by
$$ \frak p = \frak a \oplus \sum_{\alpha \in \Sigma} \lbrace w - \theta(w) 
          \vert w \in \frak{g}^{\alpha} \rbrace.$$
In case of a Riemannian symmetric space of non-compact type,
the eigenvalues $- \alpha(H)^2$ of $- (\mathrm{ad} H)^2$ 
are non-positive.

Since the Riemann curvature tensor is parallel,
the solutions of the Jacobi equation
$$ Y^{\prime \prime} = - R(Y, \gamma^{\prime})\gamma^{\prime}$$
are given by $t \mapsto f_j(t) v_j(t)$,
where the $f_j$ are functions
and the $v_j$ are the parallelly transported vector fields
along $\gamma$ such that $v_j(0)$ is an eigenvector
of the Jacobi operator. For an eigenvector
$v_j = w_j - \theta(w_j) \in \frak p$
of $- \mathrm{ad}(H)^2$, where $w_j \in \frak{g}_{\alpha_j}$,
with eigenvalue $\alpha_j(H) < 0$, one obtains
the fundamental solutions $\cosh(\alpha_j(H)t)$
and $\sinh(\alpha_j(H)t)$ for $f_j$.
It is understood that a root $\alpha$ can appear repeatedly here.
For $v_j$ with eigenvalue $0$, the fundamental solutions
are obviously $1$ and $t$.
An eigenbasis $v_1, \ldots, v_n$ in $\frak p$
of $- (\mathrm{ad}(H))^2$ determines a basis
$$ \lbrace Y_j^{\mathrm{hor}} \cosh(\alpha_j(H)t) v_j(t), 
           Y_j^{\mathrm{ver}} \sinh(\alpha_j(H)t) v_j(t),
             j = 1, \ldots, n \rbrace$$
for Jacobi fields along $\gamma$.
These have the properties $Y_j^{\mathrm{ver}}(0) = 0$
and $\nabla Y_j^{\mathrm{hor}}(0) = 0$,
$Y_j^{\mathrm{hor}}(0) = \nabla Y_j^{\mathrm{ver}}(0) = v_j$.

Jacobi fields along $\gamma$ are in 1-1 correspondence 
to vector fields along $d \gamma$ invariant under the
geodesic flow and the fibre multiplication
$ N_{\sigma}: TM \rightarrow TM, v \mapsto \sigma \cdot v$
for $\sigma \in \R$.
The fields $\xi_j$ corresponding to 
$Y_j^{\mathrm{hor}}$ and $\eta_j$
corresponding to $Y_j^{\mathrm{ver}}$
are a frame of $T(TM)$ along 
$d \gamma \vert_{\C \setminus \R}$.
In \cite{sz}, it is shown that the almost complex
tensor for the adapted complex structure
along $d \gamma$ w. r. t. the frame
$\xi_1, \eta_1, \xi_2, \eta_2, \ldots, \xi_n, \eta_n$
is given by $n$ blocks. 
A block for $v_j$ with non-zero eigenvalue 
is given by
$$
J_j = \left( \begin{array}{cc}
                - \frac{ \mathrm{Re} g_j}{ \mathrm{Im} g_j} & - 
                          \mathrm{Im} g_j - 
                   \frac{ (\mathrm{Re} g_j)^2 }{\mathrm{Im} g_j} \\
  \frac{1}{\mathrm{Im} g_j} & \frac{ \mathrm{Re} g_j}{ \mathrm{Im} g_j},
       \end{array} \right) $$
where $g_j(t + is) = \frac{1}{\alpha_j(H)} \tanh( (t+is) \alpha_j(H))$.
Inspection of the $J_j$ shows that there are
poles at $sH$ if and only if $\alpha_j(sH) \in \frac{\pi}{2} \Z$
for some $\alpha_j$. The blocks for zero-eigenvalues
do not contribute poles.
Note that the poles at $s = 0$ are due to the fact that
the $\xi_1, \eta_1, \xi_2, \eta_2, \ldots, \xi_n, \eta_n$
are not a frame there.
This confirms the result.



\section{Maximal Grauert domain of a Hermitian symmetric space of noncompact type}
Let $M=G/K$ be Hermitian symmetric of non-compact type. Thus $M$ can
be thought of as a bounded domain in some $\Bbb C^n$. 
More precisely, there is an open embedding of $M=G/K$
in its compact dual (Hermitian) symmetric space $U/K$. 
Let
$\overline{M}$ denote $M$ with the opposite complex structure. Then
$M$ embeds diagonally in $M\times \overline{M}$ and is the fixed
point set of the antiholomorphic involution $(x,y)\mapsto (y,x)$.
It is the aim of this section to show that $M\times \overline{M}$
is actually the maximal Grauert domain corresponding to $M$.
We will follow very closely Chapter VIII from the book \cite{helg}.

We use similar notation to the introduction, but now let
$\frak g^{\Bbb C}$, $\frak t^{\Bbb C}$ and $\frak p^{\Bbb C}$ denote the complexifications
of $\frak g$, $\frak t$ and $\frak p$ respectively. Also, let
$\frak h$ be a maximal abelian subalgebra of $\frak t$, then the
complexification $\frak h^{\Bbb C}$ is a Cartan subalgebra of $\frak g^{\Bbb C}$.
Let $\Delta$ be the non-zero roots of $\frak g^{\Bbb C}$ with respect to
$\frak h^{\Bbb C}$ and $\{H_{\alpha}|\alpha\in \Delta\}$ be the elements of
$\frak h^{\Bbb C}$ uniquely defined so that $\alpha(H)=B(H,H_{\alpha})$
for all $H\in \frak h^{\Bbb C}$, where $B$ is the Killing form. We let
${\frak g^{\Bbb C \alpha}}$ be the corresponding root subspaces.

Now, there exists a collection 
$\Gamma \subset \Delta$ of strongly
orthogonal roots $\Gamma = \lbrace \gamma_1, \ldots \gamma_r \rbrace$ 
and corresponding vectors $X_{\gamma}\in \frak g^{\Bbb C \gamma}$
such that 
$$\frak a=\sum_{\gamma\in\Gamma} \Bbb R (X_{\gamma}+X_{-\gamma})$$
is a maximal abelian subspace of $\frak p$.
The restricted roots $\Sigma$ in the decomposition
(\ref{root-decomposition}) are given by (case $C_r$): $$\Sigma =
\lbrace \pm \frac{\gamma_i + \gamma_j}{2}, 1 \leq j \leq i \leq r \rbrace \cup \lbrace \pm  \frac{\gamma_i - \gamma_j}{2},
                             1 \leq j < i \leq r \rbrace,$$
\noindent or (case $BC_r$): $$\Sigma = \lbrace \pm \frac{\gamma_i +
\gamma_j}{2}, 1 \leq j \leq i \leq r \rbrace \cup \lbrace \pm
\frac{\gamma_i - \gamma_j}{2}, 1 \leq j < i \leq r \rbrace \cup
\lbrace \pm \frac{\gamma_i}{2} \rbrace.$$

For this see \cite{amrt} or \cite{w}. 
This description yields the fact that the strongly orthogonal roots 
determine the domain $\omega \subset \frak a$ in the hermitian
symmetric case:
$$ \omega = \lbrace H \in \frak a: \vert \gamma(H) \vert < \frac{\pi}{2},
                    \gamma \in \Gamma \rbrace.$$
To see the inclusion  
$\lbrace H \in \frak{a}: \vert \gamma_i(H) \vert
< \frac{\pi}{2} \forall i= 1, \ldots, r \rbrace 
\subset  \lbrace H \in \frak{a}: \vert \alpha(H) \vert 
< \frac{\pi}{2} \forall \alpha \in \Psi_{\R} \rbrace$,
take an arbitrary 
$H \in \lbrace H \in \frak{a}: \vert \gamma_i(H) \vert
< \frac{\pi}{2} \forall i= 1, \ldots, r \rbrace $.
Then, if $\alpha = \frac{1}{2}( \pm \gamma_i \pm \gamma_j)$,
$\vert \alpha(H) \vert \leq \frac{1}{2}( \vert \gamma_i(H) \vert +
   \vert \gamma_j(H) \vert ) \leq \frac{1}{2} \cdot 2 \cdot \frac{\pi}{2}
                               = \frac{\pi}{2}.$
And if $\alpha = \pm \frac{1}{2} \gamma_i$,
$ \vert \alpha(H) \vert \leq \frac{1}{2} \vert \gamma_i(H) \vert
         \leq \frac{\pi}{4}$.
The opposite inclusion is obvious since $\Gamma \subset \Sigma$.

Since $\gamma_i( X_{\gamma_j} + X_{- \gamma_j}) = 2 \delta_{ij}$,
one obtains the descriptions
\begin{equation}
\omega = \lbrace H = \sum_{j=1}^r t_j \cdot
           ( X_{\gamma_j} + X_{- \gamma_j}): t_j \in 
                               (- \frac{\pi}{4} \frac{\pi}{4}) \rbrace
\end{equation}
and
\begin{equation}
\max_{\alpha\in\Sigma}|\alpha(X)|=2\max_{\gamma\in\Gamma}|t_{\gamma}|.
\end{equation}

We now use the representation from \cite{helg} of $M$ as a bounded
symmetric domain $D$. The geodesic $\sigma$ with $\sigma'(0)=X$ can
be complexified to a proper map $\sigma:\{|Im z|<r\}\to D\}$, where
$r=\frac{\pi}{4\max_{\gamma\in\Gamma} |t_{\gamma}|}$.
This follows from Corollary $7.18$.
In other words, the complexified geodesic exists on $\{|Im z|<1\}$
provided that $\max_{\alpha\in\Sigma}|\alpha(X)|=\frac{\pi}{2}$.

We can now construct the required biholomorphism from the maximal
Grau\-ert domain 
$G\times_K (Ad(K)\omega)\subset TM$ to 
$M\times\overline{M}.$

First define the map as follows. Map each geodesic in $M\subset TM$ into
$M\times\overline{M}$ diagonally. This map extends analytically to the
complexified geodesic. The above result shows that the complexified
geodesic maps into $M\times\overline{M}$ as a properly embedded disk.
Doing this on each geodesic gives us our map.
Considering the open embeddings of $M$ and $\overline{M}$
in $U/K$ and $\overline{U/K}$ respectively,
recall that the complexification $G^{\C} = U^{\C}$ acts 
transitively and holomorphically
with parabolic isotropy groups $P$ and $\tilde{P}$ 
respectively.
So our map is given by
$$
\begin{array}{rcl}
   \Omega_{\mathrm{AG}} & \rightarrow & M \times \overline{M} \\
       g \exp(i v) K^{\C} & \mapsto & (g \exp(iv) P, g \exp(iv) \tilde{P}).
\end{array} $$
In particular, it is holomorphic.

We check that the map is injective. If it were not injective, this
would correspond
to two distinct geodesics $\sigma_1$ and $\sigma_2$ in
$M\hookrightarrow M\times\overline{M}$ whose complexifications intersect
off the diagonal, say at a point $p$. Using the antiholomorphic involution,
we see that the complexified geodesics actually intersect at a second
point $p'$ off the diagonal also. We note that in a Hermitian symmetric
space a complexified geodesic corresponding to an initial geodesic $\sigma$
with $\sigma'(0)=X$ lies inside the image of the exponential map at
$\sigma(0)$ applied to the plane spanned by $X$ and $iX$. Let 
$q_1=\sigma_1(0)$ and $q_2=\sigma_2(0)$. Then there is a geodesic 
$\delta$ in
$M\times\overline{M}$ from $q_1$ to $p$ and we must have
$exp_p(\delta'(p),i\delta'(p))=exp_{q_1}(\sigma_1'(0),i\sigma_1'(0))$
as subsets of $M\times\overline{M}$.
Similarly we can find a $Y\in T_p(M\times\overline{M})$ such that
$exp_p(Y,iY)=exp_{q_2}(\sigma_2'(0),i\sigma_2'(0))$.

If the complex subspace spanned by $Y$ in $T_p(M\times\overline{M})$ is equal
to that spanned by $\delta'(p)$ then the original geodesics must coincide,
contrary to our hypothesis. Otherwise, since the complexified geodesics
also intersect at $p'$, we see that $exp_p$ cannot be a diffeomorphism,
which is also a contradiction. 
Hence, this map is injective and therefore biholomorphic
onto its image.

For surjectivity, let any point $(p,q) \in M \times \overline{M}$
be given. There is a point $(m,m)$ in the diagonal with minimal
distance $d$ to $(p,q)$ and the length minimizing geodesic 
$\sigma: \lbrack 0, d \rbrack \rightarrow M \times \overline{M}$
with $\sigma(0) = (m,m)$ and $\sigma(d) = (p,q)$
is orthogonal to the diagonal in $(m,m)$.
It follows that the geodesic $\gamma: \R \rightarrow  M \times \overline{M}$
with $\gamma(t) = \mathrm{Exp}(J \sigma^{\prime}(0))$,
where $J$ is the almost complex structure of the
K\"ahler manifold $ M \times \overline{M}$
is in fact a geodesic in the totally geodesic diagonal.
Thus $(p,q)$ is in the image of 
$d \gamma$
under embedding of $\Omega_{\mathrm{AG}}$ into
$M \times \overline{M}$.

Hence the map $G\times_K (Ad(K)\omega) \to M\times\overline{M}$ 
is a biholomorphism. 

\section{Proof of Theorem $7$}

Suppose that $N$ is an Hermitian space of the noncompact type with
an antiholomorphic involution $\sigma$ having fixed point set $M$.
Then $M$ is totally geodesic in $N$ and is itself a Riemannian
symmetric space with the restricted metric. In fact, if
$N=G'/K'$ then $M=G/K$ where $G$ is the centralizer of $\sigma$ in
$G'$ and $K$ is the centralizer of $\sigma$ in $K'$.

We will now use the notation $\frak g'$, $\frak a'$, $\frak g$, $\frak a$
in the introduction for the Lie algebras and maximal abelian subalgebras
of $N$ and $M$ respectively. We let $\Sigma$ be a root system for $\frak a$ and
$\Sigma'$ a root system for $\frak a'$.

First suppose 
$$\max_{\alpha\in\Sigma}|\alpha (H)|=\max_{\alpha\in\Sigma'}|\alpha (H)| \forall H\in \frak a.$$
In this case we want to construct a biholomorphism from $N$ to the maximal
Grauert tube over $M$. To do this, let $\gamma \subset M \subset N$ be
a goedesic. Assume that $\gamma '(0)=X\in \frak a$ and
$\max_{\alpha\in\Sigma'}|\alpha (X)|=\frac{\pi}{2}$. Then the
calculation in the previous section shows that $\gamma$ can be complexified
in $N$ to a proper map on $\{|Im z|<1\}\subset \Bbb C$. Now since
$\max_{\alpha\in\Sigma}|\alpha (X)|=\frac{\pi}{2}$ also, identifying
$\gamma$ with the same geodesic in $M\subset TM$, we can extend this
identification to an analytic isomorphism between the complexified
geodesic in $N$ and the differential of $\gamma$ in the maximal
Grauert domain.

Doing this on every geodesic, the same reasoning as in the previous
section shows that we get a well-defined map from $N$ to the
maximal Grauert domain. It is a biholomorphism since it is
analytic on each geodesic.

Finally, we notice that any biholomorphism from the maximal Grauert
domain to $N$ must be exactly of this form. Namely $M\subset TM$
must map to the fixed point set of an involution and, since
$G$ will push forward to a subgroup of $G'$, geodesics must
map to geodesics. Hence the condition in Theorem $7$ is both
necessary and sufficient.

{\bf Remark}

If the condition of Theorem $7$  is not satisfied, the map constructed
above from $N$ to $TM$ is still a biholomorphism onto its image. In
this case $N$ embeds as a subdomain of the maximal Grauert domain.
Examples of this phenomenon are the Grauert domains over the
symmetric spaces $SO(2,1)/SO(2)$ and $SO(3,2)/SO(3)\times SO(2)$.
Here
$$SO(2,1)/SO(2)\subset SU(2,1)/S(U_2 \times U_1) \subset SO(2,1)/SO(2) \times SO(2,1)/SO(2)$$
and
$$SO(3,2)/SO(3)\times SO(2)\subset SU(3,2)/S(U_3 \times U_2)$$ $$\subset SO(3,2)/SO(3)\times SO(2) \times SO(3,2)/SO(3)\times SO(2).$$
The first case here is a rank one example and so the maximal Grauert
domain coincides with the maximal Grauert tube. Normalizing the curvature
of the hyperbolic plane $H=SO(2,1)/SO(2)$ to be $-1$, we have that
$r_{max}=\frac{\pi}{2}$ and the tubes $T^{\frac{\pi}{4}}H=SU(2,1)/S(U_2 \times U_1)=B\subset \Bbb C^2$ and $T^{\frac{\pi}{2}}H=SO(2,1)/SO(2)\times SO(2,1)/SO(2)=D\times D\subset \Bbb C^2$.

\section{Proof of Proposition $8$}

We use the same notation as in the previous section.

We first show that we can assume the Hermitian symmetric space $N$ is
irreducible. Let $N = N_1 \times \ldots \times N_n$ be the
decomposition of $N$ into irreducible factors. Since the conjugation
$\sigma$ induces an isometry of $N$, it must permute the factors $N_i$
isometrically. Hence, we can number the factors so that
$\sigma(N_{2j-1}) = N_{2j}, j = 1, \ldots, k,$ and $\sigma(N_j) = N_j,
j = 2k+1, \ldots, n.$ If $M_j$ denotes the fixed point set of $\sigma$
restricted to $M_{2j-1} \times M_{2j}, j = 1, \ldots, k,$ or $M_j, j =
2k+1, \ldots, n,$ it follows that $M$ is isometric to the product of
all these $M_j$, and that $M_j$ is Hermitian and diagonal in $N_{2j-1}
\times N_{2j}, j = 1, \ldots, k,$ where $N_{2j-1}$ and $N_{2j}$ are
the same symmetric space with opposite complex structures. Note that
$N$ is a maximal Grauert domain if and only if each of the $N_j, j =
2k+1, \ldots, n$ are maximal for the corresponding $M_j$'s, since we
have already shown in section 3 that $N_{2j-1} \subset N_{2j-1}
\times N_{2j}$ as above is maximal. Thus, we may assume $N$ is
irreducible.

Now suppose that $r = rank(M) = rank(N)$ but $N$ is the
maximal Grauert domain for $M$.

The Cartan decomposition $\frak{g}'=\frak{k}'+\frak{p}'$ for $N$
with respect to the Cartan involution $\theta$
yields the decomposition $\frak{g}=\frak{k}+\frak{p}$ for $M$
where $\frak{k}=\frak{k}'\cap \frak{g}$ and $\frak{p}=\frak{p}'\cap \frak{g}$.

Choose a maximal abelian subspace $\frak{a} \subset \frak{p}$. Without
loss of generality, we can assume that $\frak{a} \subset \frak{a}'$,
a maximal abelian subalgebra of $\frak{p}'$.

Since $N$ is irreducible, there are strongly orthogonal roots
$\gamma_1,...,\gamma_r$ for $\frak{g}'$ with respect to
$\frak{a}'$. In the root decomposition

\begin{equation}
\frak{g}'=Z(\frak{a}')+\sum_{\alpha \in \Sigma'}\frak{g}'^{\alpha}
\end{equation}
it is known that $\dim(\frak{g}'^{\gamma_i})=1$ for all $i$, see for
instance \cite{amrt}, section $2.3$.

Now, associated to $\gamma_i$ there is a subalgebra $\frak{sl}_2(\Bbb
R)\subset \frak{g}'$ spanned by $\frak{g}'^{\gamma_i}$,
$\frak{g}'^{-\gamma_i}=\theta(\frak{g}'^{\gamma_i})$ and some
$x_i \in \frak{a}'$, see again \cite{amrt}. Upon exponentiation, this
$\frak{sl}_2(\Bbb R)$ generates a complex disk in $N$.

The root space decomposition of Hermitian symmetric spaces discussed
in section $3$
implies that

\begin{equation}
|\gamma_i(x_i)|>|\alpha(x_i)| \forall \alpha \in \Sigma' \setminus \{\pm \gamma_i\}.
\end{equation}

If the maximality condition holds, then there exists a $\beta \in \Sigma$
such that $\beta (x_i)=\gamma_i(x_i)$ and hence $\beta_i=\gamma_i$.
Since $\dim(\frak{g}'^{\gamma_i})=1$, we have that
$\frak{g}'^{\gamma_i}\subset \frak{g}$ and also
$\frak{g}'^{-\gamma_i}\subset \frak{g}$. Therefore the $\frak{sl}_2(\Bbb R)$
associated to $\gamma_i$ is a subalgebra of $\frak{g}$. This contradicts
the fact that $M$ is a totally-real submanifold of $N$.

Now suppose that $s=rank(N)>rank(M)=r$.
Let $J$ be the complex structure on $\frak{p}'=T_0(N)$. 
Then $\frak{a}'=\frak{a}+\frak{a}''$ where $\frak{a}''$ is contained
in the $(-1)$-eigenspace of the complex conjugation $\sigma$ on $N$ and
$J\frak{a}''\subset \frak{p}$. Here again $\frak{a}'$ is a maximal
abelian subspace of $\frak{p}'$ and $\frak{a}$ is maximal abelian
in $\frak{p}$.

The polydisk theorem (see for instance \cite{amrt} again) tells us that
$\exp(\frak{a}')$ can be complexified to a totally geodesic, embedded
polydisk $\Delta ^s \subset N$. The involution $\sigma$ restricts to
an involution on $\Delta ^s$ whose fixed-point set $P$ is totally
geodesic in $\Delta ^s$ and hence $M$. Clearly
$T_0(P)=\frak{a}+J\frak{a}''$ and $P$ has rank $r$ and dimension $s$,
since $\frak{a} \cap J \frak{a}'' = 0$: otherwise the bounded
symmetric domain would contain a $J$-invariant complete 2-dimensional
flat, which would be a copy of ${\Bbb C}$, contradicting
Liouville's theorem.

Now, an involution on a polydisk $\Delta ^s$ either preserves a factor
or preserves a pair of factors $\Delta ^2$, permuting the pair. In the
second case, up to choosing coordinates $(z,w)\in \Delta ^2$ we may
take $\sigma (z,w)=(\overline{w}, \overline{z})$ and the fixed-point
set is a copy of the hyperbolic plane $H$. Thus $P=H^p \times \Bbb R^q$.
The $\Bbb R$-factors come from factors of $\Delta ^s$ fixed by
$\sigma$ and so we may assume that they are tangent to some
$x_i \in \frak{a}$, as above, corresponding to a strongly orthogonal
root. Comparing ranks and dimensions we find that $p+q=r$ and
$2p+q=s$ respectively. Thus $p=s-r$ and $q=2s-r$.

Suppose that $s<2r$. Then $P$ contains at least one $\Bbb R$-factor
and so we may assume that there exists an $x_i \in T_0(P)\subset
T_0(M)=\frak{p}$. If $N$ were the maximal Grauert domain then we could
find an $X\in \frak{g}$ such that $[X,x_i]=\gamma_i (x_i)X$. We claim
that $X\in \frak{g}'^{\gamma_i}$.  But by the direct sum decomposition
$(4)$ above we have
$$X=z+ \sum_{\alpha \in \Sigma'}Y_{\alpha}$$
where $Y_{\alpha}\in \frak{g}'^{\alpha}$. Hence
$$\gamma_i(x_i)z+ \sum_{\alpha \in \Sigma'}(\gamma_i(x_i)-\alpha(x_i))Y_{\alpha}=0.$$
By the inequality $(5)$ this gives us that $z=0$ and $Y_{\alpha}=0$ for
$\alpha \ne \gamma_i$, justifying our claim. Now we can argue as above
to produce a holomorphic disk in $M$ and derive a contradiction as before.

Finally suppose that $s=2r$. (Certainly $s \le 2r$ as $J\frak{a}''$ is
abelian in $\frak{p}$, and so is of dimension $\leq r$
.) In this case we would like to show
that $N$ is the maximal Grauert domain. To do this, given a vector
$X\in \frak{p}$, we need to find a Jacobi field along the
corresponding geodesic which when complexified has a pole at the
boundary of the complexified geodesic (see section $2$). We may assume
that $X\in \frak{a}$ and clearly it would suffice to find such a
Jacobi field tangent to $P$. Since $s=2r$, by the argument in the
preceding paragraph, we have that $P=H^r$. Each copy of $H$ in this
factorization corresponds to two strongly orthogonal roots, say,
$\gamma_1, \gamma_2$ permuted up to sign by $\sigma$. The diagonal
copy of $Sl_2(\mathbf{R})$ in the product of the $Sl_2(\mathbf{R})$'s
corresponding to $\gamma_1, \gamma_2$ acts transitively on $H$, and
the product group acts transitively on the corresponding bidisk
$\Delta^2$ in our maximal $\Delta^s$. As each such $H$ is Hermitian,
from section $3$ we know that the correponding bidisk $\Delta ^2$ is
its maximal Grauert domain, as $\Delta^s$ is of $P = H^r$. Finally,
since $P$ is totally geodesic both in $M$ and in $\Delta^s$ which in
turn are totally geodesic in $N$, a Jacobi field (in $P$ or in $M$,
equivalently) does exist tangent to $P$ as required.

\section{Proof of Theorem $6$}

Let $M=G/K$ be an irreducible Riemannian symmetric space of noncompact type
and $N$ an associated Grauert domain. It is the aim of this section
to show that $N$ is either rigid or Hermitian symmetric.

We will use the fact that as $N$ admits a bounded strictly plurisubharmonic
function, namely the length-squared function $\rho^2$, it is a hyperbolic
complex manifold, see \cite{sib}, Theorem $3$. Let $d$ denote the
metric space structure on $N$ given by integrating the infinitesimal
Kobayashi metric. 

First note that any biholomorphism of $N$ which preserves $M$ must be
the differential of an isometry of $M$. Such a result was proven in
\cite{km}, section $6$, in the context of Grauert tubes of less than
maximal radius. The proof, however, extends to our more general case.
We remark here that the tubes of less than maximal radius are all
complete hyperbolic. We do not know how to prove this for general
Grauert domains (except as a consequence of some of them being
Hermitian symmetric). The proof of \cite{km} as written does use
tautness of the tubes, but in fact hyperbolicity is enough, see
\cite{kob}, Chapter V, Theorem $3.3$.

Suppose that $N$ is not rigid. Then we can say that there exists a
biholomorphism of $N$ which moves $M$ off itself. According to
\cite{kob}, Chapter V, Theorem $2.1$, the biholomorphism group
$Aut_{\Bbb C}(N)$ is a Lie group with compact isotropy groups.
As the orbit of a point in $N\setminus M$ under the action of $G$
has higher dimension than $M$, we deduce that the identity
component $Aut_{\Bbb C}^0(N)$ must itself move $M$ off itself.
We think of elements of the Lie algebra of $Aut_{\Bbb C}^0(N)$ as
vector fields on $N$. Then given a base point $p\in M$, there exists
a vector field in the Lie algebra transverse to $M$ at $p$. We now
use a result of Cartan and the fact that $M$ is irreducible to
find that in fact vector fields in the Lie algebra must span
$T_p(TM)$, because of course the Lie algebra is invariant under
the action of $Ad(K)$.

As a consequence of this, there exist elements of $Aut_{\Bbb C}^0(N)$
taking $p$ to all points in a sufficiently small neighbourhood $U$
of $p$ in $N$. We claim this implies that $(N,d)$ must be complete.

{\bf Proof of claim}

Again we argue by contradiction. If $(N,d)$ were not complete, there
exists a non-empty collection $\cal C$ of Cauchy sequences $\{x_i\}$
with no convergent subsequences.

Hence the function $f:N\to [0,\infty)$ given by
$$f(q)=\inf_{\cal C}\sup_i d(q,x_i)$$
is in fact finite-valued.

As $d$ is invariant under $Aut_{\Bbb C}(N)$, so is $f$ and hence
$f$ is constant, say equal to $D$, in $U$. Choose $\epsilon >0$ such that
$S=\{q\in N|d(p,q)=\epsilon\}$ is a compact subset of $U$.

Now let $\{x_i\}\in \cal C$ be such that $\sup_i d(p,x_i)<D+\frac{\epsilon}{2}$.

Every path from $p$ to $x_i$ must pass through a point in $S$. Therefore
we can find a sequence $q_i \in S$ with $d(q_i,x_i)<D-\frac{\epsilon}{2}$.

As $S$ is compact, after taking subsequences we may assume that $q_i$
converges to a point $q\in S$.

Then replacing the subsequence $\{x_i\}$ by $\{x_i\}_{i\ge N}$ for some
large $N$, we have
$$\sup_i d(q,x_i) \le \sup_i(d(q,q_i)+d(q_i,x_i)) < D-\frac{\epsilon}{4}.$$

But this contradicts $f(q)=D$, completing the proof of the claim.

We already know that the set of points in $N$ which can be mapped to
$p$ by an element of $Aut_{\Bbb C}(N)$ is open, but now since $N$ is
complete hyperbolic and in particular taut, this set is also closed.
Hence $Aut_{\Bbb C}(N)$ acts transitively.
Also, the geodesic symmetry at $p$ of $M$ extends to an involution
in $Aut_{\Bbb C}(N)$, and the fixed point set of the action by
conjugation is exactly the (compact) isotropy group of $p$. Thus $N$
is an Hermitian symmetric space as required.

\section{Direct Proof that the Akhiezer-Gindikin Domain is Stein.}

\begin{theorem}
Let $\frak{g} = \frak{k} \oplus \frak{p}$ be the Cartan
decomposition, $\frak{a} \subset \frak{p}$ maximal abelian,
$W$ the Weyl group for the adjoint action of $K$
on $\frak{p}$
and $\omega \subset \frak{a}$ the set describing the
Akhiezer-Gindikin domain 
$$\Omega_{\mathrm{AG}} = \lbrace g \exp(i \xi) K^{\Bbb C}, g \in G, 
                            \xi \in \omega \rbrace.$$
Let $u$ be a strictly convex, $W$-invariant, smooth function on $\omega$
and $\tilde{u}$ its $G$-invariant extension on $\Omega_{AG}$.
Then $\tilde{u}$ is strictly plurisubharmonic.
\end{theorem}

{\bf Proof}

Let us set the notation. 
The symbols $\frak{g}, \frak{k}, \frak{p}, \frak{a},
\theta$ denote the Lie algebras, Cartan involution, eigenspaces,
etc. as above, associated to our symmetric space $M$. Let $\frak{g}^{\Bbb C} =
\frak{g} \otimes {\Bbb C}, \frak{k}^{\Bbb C}, \theta,$ etc., be
the corresponding complexified objects. Let $G, K, G^{\Bbb C}, K^{\Bbb C}$,
etc., be the corresponding groups and complex groups. Let $W$ be the
Weyl group of $\frak{g}$ with respect to $K, \frak{a}$. Let $X = G^{\Bbb
C}/K^{\Bbb C}$ be the affine complexification of $M$. Let $\frak{u}$
be the compact twin of $\frak{g}$, where $\frak{u} = \frak{k} \oplus i
\frak{p}$, where multiplication by $i$ is meant in $\frak{g}^{\Bbb
C}$, and let $U \subset G^{\Bbb C}$ be the corresponding
subgroup, the maximal compact subgroup of $G^{\Bbb C}$. Let $\xi_1, \xi_2,...,
\xi_n$ be an orthonormal basis of $\frak{p}$. 
There is an orthogonal direct sum decomposition $\frak{p}=\frak{a}+\frak{p}'$
where $\frak{p}'$ is the orthogonal complement of $\frak{a}\subset \frak{p}$.
Let $M^c = U/K
\subset X$ be the compact twin symmetric space of our original $M$.
If $\sigma$ denotes the conjugation of $X$ fixing $M$, let $\tau =
\theta \circ \sigma = \sigma \circ \theta$ be the conjugation for
$M^c$. Let $\Omega_{AG}$ be the Akhiezer-Gindikin domain, which is the
image $G \cdot \exp(i Ad(K) \omega) \cdot K^{\Bbb C}$ in 
$G^{\Bbb C}/K^{\Bbb C}$. We write here $\Omega_{AG}$
for $G \times_K Ad(K) \omega$, too. 
Here the $\exp$ can be
taken either as group exponential in $G^{\Bbb C}$, or as the
geodesic exponential in the compact twin $M^c$. For any vector $\xi
\in \frak{g} \otimes {\Bbb C}$, let $\tilde{\xi}$, or sometimes $[\xi]^{\tilde{}}\;$, denote the
corresponding vectorfield given by
$\tilde{\xi}(p) = \frac{\mathrm{d}}{\mathrm{d} t} \exp(t \xi).p\mid_{t=0}$ 
on either $X$, or equivalently on $Ad(K)\omega \subset i\frak{p}$, 
for fields tangent to $M^c$ within
$\Omega_{AG}$. Let $u$ be a smooth, strictly convex, $W$-invariant
function on $\frak{a}$, and let $\tilde{u}$ be the corresponding
function on $\Omega_{AG}$. We wish to show that $\tilde{u}$ is strictly
plurisubharmonic. Fix $\xi_0 \in Ad(K) \omega$; we want to show 
$i\partial \bar{\partial}(\tilde{u}) > 0$ at $x_0 = \exp(i\xi_0) \cdot K \in
M^c$. Without loss of generality, we can assume $\xi_0 \in \omega$,
in particular, $\xi_o \in \frak{a}$.

\vskip 3mm

The only way we use that we are inside $\Omega_{AG}$ at $x_0$ is in
the following observation: $\tilde{\xi}_1,...,\tilde{\xi}_n$ are
linearly independent at $x_0$, and span a totally real subspace of the
tangent space of $X$ at $x_0$. This follows immediately from the
(equivariant) bundle structure of the Grauert domain construction.

\vskip 3mm

Since the group $G^{\Bbb C}$ acts holomorphically on $X$, we have that
$Z_k = \tilde{\xi}_k^{1,0}$ is a holomorphic vector field, and
$\bar{Z}_j$ is anti-holomorphic. As a result, $[Z_l, \bar{Z}_k] = 0$,
for all $k,l = 1,...,n.$ Since the span of the $\tilde{\xi}_k$ is
totally real, the complex fields $Z_1,..., Z_n$ are ${\Bbb
C}$-linearly independent at $x_0$, and span the complex tangent space
at that point. Thus, it suffices to show the matrix $i \partial
\bar{\partial}(\tilde{u})(Z_i, \bar{Z}_j)$ is positive definite. This
reduces immediately to the matrix $Z_i \cdot \bar{Z}_j (\tilde{u})$.

\vskip 3mm

Write $Z_k$ as $\frac{1}{2}(\tilde{\xi}_k - iJ\tilde{\xi}_k)$, and
therefore, $$Z_k \bar{Z}_l(\tilde{u}) = \frac{1}{4}(\tilde{\xi}_k -
iJ\tilde{\xi}_k)(\tilde{\xi}_l + iJ\tilde{\xi}_l)(\tilde{u})$$ $$=
\frac{1}{4}(i\tilde{\xi}_k J\tilde{\xi}_l + J\tilde{\xi}_k
J\tilde{\xi}_l)(\tilde{u}),$$ since $\tilde{\xi}_l(\tilde{u}) \equiv
0$, for all $l$. Note also that $J\tilde{\xi}_k = \widetilde{(i{\xi}_k)}$,
for all $k$. Finally, note that $$\tilde{\xi}_k
J\tilde{\xi}_l(\tilde{u}) = [\tilde{\xi}_k, J\tilde{\xi}_l](\tilde{u})
 = (i[\xi_k,\xi_l])^{\tilde{}}\;(\tilde{u})= J ([\xi_k, \xi_l])^{\tilde{}}\;(\tilde{u}),$$ since the flow of $\tilde{\xi}_k$ is
holomorphic. Now $[\xi_k, \xi_l] = \eta \in \frak{k}$, and as
$\frak{k}$ acts on $M^c$ at $x_0$, we have the basic observation that
$${\mbox{the tangent space to the orbit of}} \; K \; {\mbox{through}}
\; x_0
\; {\mbox{is spanned by vectors}} \; \widetilde{i \xi}, \xi \in \frak{p}'.$$

(Notice that these two subspaces of tangent vectors are equal at
points $x_0 = \exp(i\xi_0) \cdot 0,$ where $\xi_0 \in \frak{p}$ is a
regular element, and hence the claim holds for any $\xi_0 \in
\frak{p}$.) Hence, at $x_0$ there is some $\xi \in \frak{p}$ such that
$$(i [\xi_k, \xi_l])^{\tilde{}}\; = J ([\xi_k,
\xi_l])^{\tilde{}}\; = J \widetilde{i \xi} = -\tilde{\xi}$$ at $x_0$. But then
this implies that 
$$ \tilde{\xi_k} J \tilde{\xi_l}( \tilde{u} ) =
     [\tilde{\xi}_k, J \tilde{\xi}_l](\tilde{u}) =
     -\tilde{\xi}(\tilde{u}) = 0,$$ 
at $x_0$, again, by $G$-invariance of
$\tilde{u}$. So, we have reduced our task to showing that the matrix
$$\widetilde{i\xi_k} \widetilde{i\xi_l}(\tilde{u})$$ is positive definite at $x_0$.

\vskip 3mm

We now write the basis $\xi_1, \xi_2,...,
\xi_n$ of $\frak{p}$ so that $\xi_i,
\xi_j, j, i = 1,..., r=\dim(\frak{a})$, denote elements in $\frak{a}$,
and $\xi_k,
\xi_l$ from among $\xi_{r+1},...,\xi_n$ are in $\frak{p}' =$
orthogonal complement of $\frak{a} \subset \frak{p}$. More precisely,
let $\Sigma$ be the set of roots for $\frak{a}$ acting on $\frak{g}$,
and let $\Sigma_{+}$ be the positive roots with respect to some
ordering. For $\alpha \in \Sigma$, let $X_{\alpha}$ be a non-zero
vector in $\frak{g}^{\alpha} = \{Y \in \frak{g} \mid [H, Y] =
\alpha(H) Y, \; {\mbox{for all}} \; H \in \frak{a}\}.$ For $\alpha \in
\Sigma, X_{\alpha} - \theta X_{\alpha} \in \frak{p}',$ and in fact, a
basis for $\frak{p}'$ is given by $X_{\alpha} - \theta X_{\alpha},
\alpha \in \Sigma_{+}.$ Note that in this notation we are counting
the $\alpha's \in \Sigma_{+}$ with {\em multiplicity}. Finally, the
$X_{\alpha} - \theta X_{\alpha}$ are orthogonal to one another, and
so, up to scale, we may take them as an orthonormal basis for
$\frak{p}'$. In other words, for $k, l = r+1,...,n,$ we can take
$\xi_k, \xi_l$ to be of the form $X_{\alpha} - \theta X_{\alpha}$, up
to a scale factor. Let $x_1,...,x_n$ denote the affine coordinates in
$\frak{p}$ corresponding to this (orthonormal) basis of $\frak{p}$.

\vskip 3mm

We have three types of terms to consider in the matrix $\widetilde{i\xi_s}
\widetilde{i\xi_t}(\tilde{u})$:

\vskip 2mm

\hspace*{0.5in} i) $\frak{a}$-terms $\widetilde{i\xi_i}
\widetilde{i\xi_j}(\tilde{u}), i, j = 1,...,r.$

\vskip 2mm

\hspace*{0.5in} ii) cross-terms $\widetilde{i\xi_i}
\widetilde{i\xi_k}(\tilde{u}), i = 1,...,r, k = r+1,...,n.$

\vskip 2mm

\hspace*{0.5in} iii) $\frak{p}'$-terms $\widetilde{i\xi_k}
\widetilde{i\xi_l}(\tilde{u}), k, l = r+1,...,n.$

\vskip 3mm
In case ii), notice that $[Z_s, \bar{Z}_t] \equiv 0$ and the calculation above imply that
$[\widetilde{i\xi_i}, \widetilde{i\xi_k}](\tilde{u}) \equiv 0$ also.

\vskip 3mm

The first two types of terms are easy to compute, and do not even
require us to know the fields $\widetilde{i\xi_s}$ very explicitly.

\vskip 3mm

{\em Case i):} In this case, the subalgebra $\frak{a}$ is abelian, and
the subspace $\exp(i\frak{a}) \cdot 0 \subset M^c$ is flat. Therefore,
$$\widetilde{i\xi_i} \widetilde{i\xi_j}(\tilde{u})(x_0) =
\frac{\partial^2 u}{\partial x_i \partial x_j} (\xi_0).$$

\vskip 3mm

{\em Case ii):} First use the basic observation above to see that, at
all regular $\xi_0 \in \frak{a}$, the tangent vector $\widetilde{i\xi_k}$
is tangent to the orbit of $K$ through $x_0 = \exp(i\xi_0) \cdot 0$,
and by $K$-invariance of $\tilde{u}$, $\widetilde{i\xi_k}(\tilde{u})(x_0)
= 0$. Since the regular $\xi_0$ are dense in $\frak{a}$, this is true
at every $x_0 \in \exp(i\frak{a})
\cdot 0$. Therefore, since every $\widetilde{i\xi_i}$ is tangent to
$\exp(i\frak{a}) \cdot 0$, we have that $\widetilde{i\xi_i}
\widetilde{i\xi_k}(\tilde{u})(x_0) = 0.$

\vskip 3mm

{\em Case iii):} Here we have to examine the differential of the
exponential map more clearly to see how the fields $\widetilde{i\xi_s}$ are
transported from $M^c$ back to $\frak{p}$ by the inverse of the
exponential map.

\vskip 3mm

We will work along constant speed geodesics $\gamma_{\xi}(t) =
\exp(it\xi) \cdot 0, \xi \in \frak{p},$ in $M^c$. Since we will be
using the Jacobi equation, it is useful to understand parallel
transport along $\gamma_{{\xi}_0}$ explicitly. Let $g_t =
\exp(it\xi_0)\in U$, note that the tangent space to $M^c$ at $g_t =
\gamma_{{\xi}_0}(t)$ is identified with ${dg_t}_{*}\frak{p} = \frak{u}
\; {\mbox{mod}} \; Ad(g_t)(\frak{k}).$ Parallel vectors along
$\gamma_{\xi_0}$ are simply those of the form ${dg_t}_{*}(\xi)$ for
fixed $\xi \in \frak{p}$. Given this identification, a Jacobi field
$Y(t)$ can be written in terms of the parallel fields as $$Y(t) =
{dg_t}_{*}(\sum_{j = 1}^r v_j(t) i\xi_j + \sum_{\alpha \in \Sigma_{+}}
v_{\alpha}(t) i(X_{\alpha} - \theta X_{\alpha})),$$ and the coefficient
functions satisfy the ordinary differential equations: $$\ddot{v}_j
\equiv 0,$$ $$\ddot{v}_{\alpha} = \{\alpha(i\xi_0)\}^2 v_{\alpha} =
-\{\alpha(\xi_0)\}^2 v_{\alpha}.$$ Taking into account the initial
conditions $v_{\alpha}(0) = 0, \dot{v}_{\alpha}(0) = 1,$ and setting
all other coefficients equal to $0$, we get that $$d\exp_{*}(i\xi_0):
i(X_{\alpha} - \theta X_{\alpha}) \longrightarrow
{dg_1}_{*}(\frac{\sin(\alpha(\xi_0))}{\alpha(\xi_0)} i(X_{\alpha} -
\theta X_{\alpha})).$$ We are assuming, provisionally, that $\xi_0$ is
regular, so $\alpha(\xi_0) \neq 0$, for all $\alpha \in \Sigma.$ If we
denote the inverse map to $\exp$ by ``$\; {\mbox{log}} \;$'', then we
have $$d\; {\mbox{log}}_{*}(x_0) \;: {dg_1}_{*}(i(X_{\alpha} - \theta
X_{\alpha})) \longrightarrow \frac{\alpha(\xi_0)}{\sin(\alpha(\xi_0))}
i(X_{\alpha} - \theta X_{\alpha}).$$ Recall that, since $0 \neq \xi_0
\in \omega_0$, we have $0 < \mid \alpha(\xi_0) \mid < \frac{\pi}{2},$
so the denominator above doesn't vanish.

\vskip 3mm

Next we have to figure out how to represent the fields
$\widetilde{i\xi_k}$ back on $\frak{p}$. Of course, these are just the
derivatives of the actions of $\exp(it(X_{\alpha} - \theta
X_{\alpha})) \in U$ on $M^c$. This is computable at $x_0 = g_1$ as
${dg_1}_{*}[ \frac{d}{dt} \mid_{t=0} g_1^{-1}\cdot \exp(it(X_{\alpha}
- \theta X_{\alpha})) \cdot g_1 \cdot 0]$. Now this curve at $0 \in
M^c$ obviously has as derivative at $t = 0$ the image of
$Ad(g_1^{-1})i(X_{\alpha} - \theta X_{\alpha}) \in \frak{u}$ modulo
$\frak{k}$. We can calculate this using $g_1 = \exp(i\xi_0) \in
U$. So, $$Ad(g_1^{-1})i(X_{\alpha} - \theta X_{\alpha}) = i e^{-i
ad(\xi_0)} (X_{\alpha} - \theta X_{\alpha})$$ $$= e^{-i \alpha(\xi_0)}
iX_{\alpha} - i e^{i \alpha(\xi_0)} \theta X_{\alpha}$$ $$ =
\cos(\alpha(\xi_0)) i(X_{\alpha} - \theta X_{\alpha}) -
\sin(\alpha(\xi_0))(X_{\alpha} + \theta X_{\alpha})$$ $$=
\cos(\alpha(\xi_0)) i(X_{\alpha} - \theta X_{\alpha}) \; {\mbox{mod}}
\; \frak{k},$$ since $X_{\alpha} + \theta X_{\alpha} \in \frak{k}$,
for every $\alpha \in \Sigma$. From this we conclude $$d\log_{*}(x_0):
[i(X_{\alpha} - \theta X_{\alpha})]^{\tilde{}}\; \longrightarrow
\frac{\alpha(\xi_0) \cos(\alpha(\xi_0))}{\sin(\alpha(\xi_0))}
i(X_{\alpha} - \theta X_{\alpha}).$$

\vskip 3mm

Now, before computing $\widetilde{i\xi_k}
\widetilde{i\xi_l}(\tilde{u})(x_0)$, we make an auxiliary computation
of $\widetilde{i\xi_k}
\widetilde{i\xi_l}(\tilde{u}_0)(x_0)$, where $u_0$ denotes the
gradient of $u$ on $\frak{p}$ at $\xi_0$, and by abuse of notation,
will also denote the function $\xi \rightarrow u_0 \cdot \xi$
determined by taking the inner product of $\xi$ with $u_0$. Similarly,
$\tilde{u}_0$ denotes the same function transported to $\exp(iAd(K)
\omega_0) \cdot 0 \subset M^c$ by the exponential function. We carry
out the computation on $\frak{p}$ near $\xi_0$, that is, we compute
$$[i(X_{\alpha} - \theta X_{\alpha})]^{\tilde{}}\;[i(X_{\beta} -
\theta X_{\beta})]^{\tilde{}}\;(\tilde{u}_0)(x_0)$$ $$ =
d\log_{*}([i(X_{\alpha} - \theta X_{\alpha})]^{\tilde{}}\;(x_0))(u_0 \cdot
d\log_{*}([i(X_{\beta} - \theta X_{\beta})]^{\tilde{}}\;)),$$ evaluated
at $\xi_0 \in \frak{p}$.

\vskip 3mm

We first want to use the basic observation above to replace the
tangent vector $[i(X_{\alpha} - \theta X_{\alpha})]^{\tilde{}}\;$ at
$x_0$ by $\tilde{\eta}$, for suitable $\eta \in \frak{k}$. We repeat
that we are assuming provisionally that $\xi_0$ is regular. To do
this, let us calculate much as before $$(X_{\alpha} + \theta
X_{\alpha})^{\tilde{}}\;(x_0) = {dg_1}_{*}[\frac{d}{dt} \mid_{t=0} g_1^{-1} \cdot
\exp(t(X_{\alpha} + \theta X_{\alpha})) \cdot g_1 \cdot 0]$$ $$=
{dg_1}_{*}(e^{-i \; {\mbox{ad}} \; (\xi_0)}(X_{\alpha} + \theta
X_{\alpha}))$$ $$ = {dg_1}_{*}(\cos(\alpha(\xi_0))(X_{\alpha} + \theta
X_{\alpha}) - i \sin(\alpha(\xi_0))(X_{\alpha} - \theta X_{\alpha}))$$
 $$= {dg_1}_{*}(-\sin(\alpha(\xi_0))i(X_{\alpha} - \theta
X_{\alpha})) \mathrm{mod} \frak k.$$ 
Taken together with what was shown above, we conclude
$$[i(X_{\alpha} - \theta X_{\alpha})]^{\tilde{}}\;(x_0) = -
\cot(\alpha(\xi_0)) (X_{\alpha} + \theta X_{\alpha})^{\tilde{}}\;.$$
Going back to our original computation, this gives us that $$[i(X_{\alpha} - \theta X_{\alpha})]^{\tilde{}}\;[i(X_{\beta} -
\theta X_{\beta})]^{\tilde{}}\;(\tilde{u}_0)(x_0)$$ $$ = - \cot(\alpha(\xi_0))
d\log_{*}([i(X_{\alpha} + \theta X_{\alpha})]^{\tilde{}}\;(x_0))(u_0 \cdot
d\log_{*}([i(X_{\beta} - \theta X_{\beta})]^{\tilde{}}\;)).$$ Because the
subgroup $K$ fixes $0 \in M^c$, the exponential map from $0$ is
equivariant with respect to the action of $K$ on $T_0(M^c) = \frak{p}$
and on $M^c$. Thus, we conclude that 
$$[i(X_{\alpha} - \theta X_{\alpha})]^{\tilde{}}\;[i(X_{\beta} -
\theta X_{\beta})]^{\tilde{}}\;(\tilde{u}_0)(x_0)$$ $$= - \cot(\alpha(\xi_0))
\frac{d}{dt} [u_0 \cdot d\log_{*}([i(X_{\beta} -
\theta X_{\beta})]^{\tilde{}}\;)(Ad(\exp(t(X_{\alpha} + \theta X_{\alpha})))\xi_0)] 
\mid_{t=0}.$$
For simplicity, set $k_t = \exp(t(X_{\alpha} + \theta X_{\alpha})) \in
U$. $K$-equivariance of $\exp$ implies $$\frac{d}{dt} [u_0 \cdot d\log_{*}([i(X_{\beta} -
\theta X_{\beta})]^{\tilde{}}\;)(Ad(k_t)\xi_0)]
\mid_{t=0}$$ $$ = u_0 \cdot [\frac{d}{dt}(Ad(k_t) \cdot d\log_{*}(x_0)\{
{dk_t^{-1}}_{*}([i(X_{\beta} - \theta X_{\beta})]^{\tilde{}}\;)\}
\mid_{t=0}]$$ $$= u_0 \cdot [\frac{d}{dt}(Ad(k_t) \cdot d\log_{*}(x_0)\{
[iAd(k_t^{-1})(X_{\beta} - \theta X_{\beta})]^{\tilde{}}\;\}
\mid_{t=0}]$$ $$= u_0 \cdot [X_{\alpha} + \theta X_{\alpha},
d\log_{*}(x_0)([i(X_{\beta} - \theta X_{\beta})]^{\tilde{}}\;)]$$ $$ - u_0
\cdot d\log_{*}(x_0)([i[X_{\alpha} + \theta X_{\alpha},
X_{\beta} - \theta X_{\beta}]]^{\tilde{}}\;).$$

\vskip 3mm

We next compute $[X_{\alpha} + \theta X_{\alpha}, i(X_{\beta} - \theta
X_{\beta})] = i([\theta X_{\alpha}, X_{\beta}] - \theta [\theta
X_{\alpha}, X_{\beta}]).$ Notice that if $\alpha \neq \beta$, then
$i([\theta X_{\alpha}, X_{\beta}] - \theta [\theta X_{\alpha},
X_{\beta}]) \in i \frak{p}'$, which is tangent to the $K$-orbit of
$\xi_0$. Hence $u_0 \cdot i([\theta X_{\alpha}, X_{\beta}] - \theta
[\theta X_{\alpha}, X_{\beta}])) = 0.$

\vskip 3mm

For the second term above, we note that $d\log_{*}(x_0)$ sends the
subspace spanned by $[i(X_{\alpha} - \theta X_{\alpha})]^{\tilde{}}\;$ to
the subspace spanned by $i(X_{\alpha} - \theta X_{\alpha})$, and
similarly the subspace spanned by $\widetilde{i\xi_j}, j = 1,...,r,$
gets sent to the subspace $i\frak{a} \subset i\frak{p}$. As a result,
$$ u_0 \cdot d\log_{*}(x_0)(-[i[X_{\alpha} + \theta
X_{\alpha}, X_{\beta} - \theta X_{\beta}]]^{\tilde{}}\;) $$ $$
 = u_0 \cdot
d\log_{*}(x_0)([i\{[\theta X_{\alpha}, X_{\beta}] 
- \theta [\theta X_{\alpha}, X_{\beta}] + [X_{\alpha}, X_{\beta}] 
- \theta([X_{\alpha}, X_{\beta}])\}]^{\tilde{}}\;) = 0,$$ for all $\alpha \neq
\beta \in \Sigma_{+}.$ 

\vskip 3mm

Finally, in the case when $\alpha = \beta$, then $i([\theta
X_{\alpha}, X_{\alpha}] - \theta [\theta X_{\alpha}, X_{\alpha}]) \in
i Z(\frak{a}) \cap \frak{a} = i\frak{a}$. 
Recalling that $d\log_{*}(x_0)$ is the ``identity'' when
restricted to vectors tangent to the flat $\exp(i\frak{a}) \cdot 0
\subset M^c$, we get $$ u_0 \cdot d\log_{*}(x_0)([i\{[\theta
X_{\alpha}, X_{\alpha}] - \theta [\theta X_{\alpha}, X_{\alpha}]\}]^{\tilde{}}\;) =
u_0 \cdot i\{[\theta X_{\alpha}, X_{\alpha}] - \theta [\theta
X_{\alpha}, X_{\alpha}]\},$$ 
where $i\{[\theta X_{\alpha}, X_{\alpha}] - \theta [\theta
X_{\alpha}, X_{\alpha}]\} \in \frak{a}$.
Putting this together with
what we have left of the first term, we get $$[i(X_{\alpha}
- \theta X_{\alpha})]^{\tilde{}}\;[i(X_{\beta} - \theta
X_{\beta})]^{\tilde{}}\;(\tilde{u}_0)(x_0) = 0, \; {\mbox{if}} \; \alpha \neq
\beta,$$ and $$[i(X_{\alpha} - \theta
X_{\alpha})]^{\tilde{}}\;[i(X_{\alpha} - \theta
X_{\alpha})]^{\tilde{}}\;(\tilde{u}_0)(x_0)$$ $$ = u_0 \cdot i([\theta X_{\alpha},
X_{\alpha}] - \theta [\theta X_{\alpha}, X_{\alpha}]) - u_0 \cdot
i([\theta X_{\alpha}, X_{\alpha}] - \theta [\theta X_{\alpha},
X_{\alpha}]) = 0,$$ for all $\alpha \in \Sigma_{+}.$

Summarizing, we have shown $$[i(X_{\alpha} - \theta
X_{\alpha})]^{\tilde{}}\;[i(X_{\beta} - \theta
X_{\beta})]^{\tilde{}}\;(\tilde{u}_0)(x_0) = 0, \; {\mbox{for all}} \;
\alpha, \beta \in \Sigma_{+}.$$ By our choice of $u_0$ as the gradient
of $u$ at $\xi_0$, we have that the function $u - u_0 $ has a critical
point at $\xi_0 \in \frak{p}$, equivalently, $\tilde{u} - \tilde{u}_0$
has a critical point at $x_0$. (Recall that, by abuse of notation, we are
denoting by $u_0$ the linear function $u_0 \cdot \xi$.) Hence the
derivatives
$$[i(X_{\alpha} - \theta X_{\alpha})]^{\tilde{}}\;[i(X_{\beta}
- \theta X_{\beta})]^{\tilde{}}\;(\tilde{u})(x_0) = [i(X_{\alpha} -
\theta X_{\alpha})]^{\tilde{}}\;[i(X_{\beta} - \theta
X_{\beta})]^{\tilde{}}\;(\tilde{u} - \tilde{u}_0)(x_0)$$ can be evaluated in terms
of the invariantly defined Hessian of $u - u_0$ at $\xi_0 \in
\frak{p}$, that is, we get $$[i(X_{\alpha} - \theta
X_{\alpha})]^{\tilde{}}\;[i(X_{\beta} - \theta X_{\beta})]^{\tilde{}}\;(\tilde{u} -
\tilde{u}_0)(x_0) $$ $$= {\mbox{Hess}}(u - u_0)(\xi_0) \;
(d\log_{*}(x_0)([i(X_{\alpha} - \theta X_{\alpha})]^{\tilde{}}\;),
d\log_{*}(x_0)([i(X_{\beta} - \theta X_{\beta})]^{\tilde{}}\;))$$ $$=
\alpha(\xi_0) \cot(\alpha(\xi_0)) \beta(\xi_0) \cot(\beta(\xi_0)) \;
{\mbox{Hess}}(u - u_0)(\xi_0) (i(X_{\alpha} - \theta X_{\alpha}),
i(X_{\beta} - \theta X_{\beta})).$$ This last matrix is obviously
positive definite when $\xi_0$ is regular, since $u$ is strictly
convex. For $\xi_0$ not regular, the result follows from the regular
case, passing to the limit $\xi_0$ from regular $\xi \in \frak{a}$,
taking into account that the function
$$\alpha(\xi) \cot(\alpha(\xi)) = \frac{\alpha(\xi)}{\sin(\alpha(\xi))} \cdot
\cos(\alpha(\xi))$$ has a finite, non-vanishing limit at $\xi_0$ as
long as $\mid \alpha(\xi_0) \mid < \frac{\pi}{2}$. This completes
the proof of the theorem.

\begin{corollary}
The domain $\Omega_{AG}$ has the Stein property. 
\end{corollary}

{\bf Proof}

Let $\Gamma$ be a discrete subgroup
of $G$ acting completely discontinuously on $M = G/K$
such that $\Gamma \backslash M$ is a compact manifold.
Such a $\Gamma$ can be found according to \cite{borel}.
Let $\tilde{u}$ be given as in the theorem
by a strictly convex, $W$-invariant function
$u$ compactly exhausting $\omega$.
For instance, $$u(\xi) :=  \sum_{\alpha \in \Sigma}
\frac{1}{(\frac{\pi}{2})^2 - \alpha(\xi)^2}$$
has these properties with respect to the symmetric space metric
on $\frak p \cong T_{eK} G/K$ restricted to $\frak a$. 
Then $\tilde{u}$ pushes down to a smooth, strictly
plurisubharmonic exhaustion function of 
$\Gamma \backslash \Omega_{\mathrm{AG}}$, which is therefore Stein.
The universal covering $\Omega_{\mathrm{AG}}$
of $\Gamma \backslash \Omega_{\mathrm{AG}}$  
is then also Stein (see \cite{stein}, p. 66 f.)
completing the proof of the corollary.

\end{document}